\documentclass[11pt]{amsart}
\usepackage{geometry}                
\usepackage{epsfig}
\usepackage{psfrag}
\usepackage{amssymb}
\usepackage{epstopdf}
\usepackage{color}
\usepackage{labelfig}
\newtheorem{theorem}{Theorem}[section]    

\theoremstyle{definition}
\newtheorem{definition}[theorem]{Definition}
\newtheorem{remark}[theorem]{Remark}  
\newtheorem{example}[theorem]{Example}    
\numberwithin{equation}{section}

\newcommand{\F}{\mathcal F_{ob} }

\newcommand{\Aut}{\rm{Aut}}

\usepackage{fancyhdr}
\begin{document}

\title{Visualizing overtwisted discs in open books}
\author{Tetsuya Ito}
\address{Research Institute for Mathematical Sciences, Kyoto university
Kyoto, 606-8502, Japan}
\email{tetitoh@kurims.kyoto-u.ac.jp} 
\urladdr{http://kurims.kyoto-u.ac.jp/~tetitoh/}
\author{Keiko Kawamuro}
\address{14 MacLean Hall\\ Department of Mathematics \\ 
The University of Iowa\\ Iowa City, IA\\ 52242-1419, USA}
\email{kawamuro@iowa.uiowa.edu}


\date{\today} 

\begin{abstract}
We give an alternative proof of a theorem of Honda-Kazez-Mati\'c that every non-right-veering open book supports an overtwisted contact structure. 
We also study two types of examples that show how overtwisted discs are embedded relative to right-veering open books. 
\end{abstract}

\maketitle

\section{Introduction}

In \cite{ik1-1}, we have introduced open book foliations and their basic machinery by using that of braid foliations \cite{BM4, BM2, BM5, BM1, BM6, BM3, bm1, bm2} and showed applications of open book foliations including a self-linking number formula of general closed braids. 
In \cite{ik2} we study the geometric structure of a 3-manifold by using open book foliations. 
In this note we study more applications of open book foliations. 
We will assume the readers are familiar with the definition and basic machinery of open book foliations in \cite{ik1-1}.

One of the features of open book foliations is that one can visualize how surfaces are embedded with respect to general open books. 
In this paper we use this feature to illustrate overtwisted discs and give constructive methods to detect overtwisted contact structures.

We first give an alternative proof of a tightness criterion theorem by Honda, Kazez and Mati\'c \cite{hkm1}: If an open book is not right-veering then it supports an overtwisted contact structure. 

The converse does not hold: 
In fact, Honda, Kazez and Mati\'c \cite{hkm1} show that if a contact structure $\xi$ is supported by a non-right veering open book $(S, \phi)$, by applying positive stabilizations to $(S, \phi)$ one can find a right-veering open book $(\widehat S, \widehat \phi)$ that supports $\xi$.   
We concretely visualize an overtwisted disc relative to the right-veering $(\widehat S, \widehat \phi)$.

Lastly, we give an infinite family of open books that are right-veering and non-destabilizable but compatible with overtwisted contact structures.
This negatively answers a question of Honda, Kazez and Mati\'c \cite{hkm1}. 
Our family generalizes the previously known examples by Lekili \cite{le} and Lisca \cite{lis}, but our proof of overtwistedness is more direct.

\section{Overtwisted discs in non-right-veering open books}

Recall that an {\em overtwisted disc} is an embedded disc in a contact 3-manifold whose boundary is a limit cycle in the characteristic foliation of the disc. In particular, every  overtwisted disc has {\em Legendrian} boundary. In the framework of open book foliations a {\em transverse overtwisted disc} plays a corresponding role:

\begin{definition}\label{def:trans-ot-disc}
\cite[Def 4.1]{ik1-1} 
Let $D \subset M_{(S, \phi)}$ be an oriented disc whose boundary is positively braided (i.e., a transverse knot) with respect to the open book $(S, \phi)$. 
If the following are satisfied $D$ is called a {\em transverse overtwisted disc}: 
\begin{enumerate}
\item $G_{--}$ is a connected tree with no fake vertices,
\item $G_{++}$ is homeomorphic to $S^1$.
\item $\F(D)$ contains no c-circles,
\end{enumerate}
(Terminologies like $G_{--}, G_{++}$, fake vertices and c-circles are defined in \S 2.1 of \cite{ik1-1}.)
\end{definition}

In \cite{ik1-1} we show that the manifold $M_{(S,\phi)}$ contains a transverse overtwisted disc if and only if the  contact manifold $(M_{(S,\phi)}, \xi_{(S, \phi)})$ contains an overtwisted disc. Hence from now on, we may {\em not}  distinguish a transverse overtwisted disc and a usual overtwisted disc, and often call a transverse overtwisted disc simply an overtwisted disc.

Next we review the notion of right-veering mapping classes then reprove Honda Kazez Mati\'c's tightness criterion in Theorem~\ref{thm:nrv-ot}.

\begin{definition}\cite{hkm1} 
Let $\gamma, \gamma'$ be oriented properly embedded arcs in the surface $S$ that  start from the same point $* \in \partial S$. 
Suppose, after some isotopy relative to the endpoints, $\gamma$ and $\gamma'$ realize the minimal geometric intersection number. 
We say that $\gamma'$ lies strictly {\em on the right side} of $\gamma$ if around the common starting point $*$ the curve $\gamma'$ strictly lies on the right side of $\gamma$. In such case we denote $\gamma > \gamma'$.
\end{definition}

\begin{definition}\cite[Definition 2.1]{hkm1}
\label{defn:right-veering}
Let $C$ be a boundary component of $S$.
We say that $\phi \in \Aut(S, \partial S)$ is {\em right-veering} with respect to $C$ if $\gamma \geq \phi(\gamma)$ holds for any isotopy classes $\gamma$ of properly embedded curves which start at a point on $C$.
We say that the diffeomorphism $\phi$ (or the open book $(S, \phi)$) 
is {\em right-veering} if $\phi$ is right-veering with respect to all the boundary components of $S$. 
In particular, the identify ${\rm{id}} \in {\Aut}(S, \partial S)$ is right-veering.
\end{definition}

The following theorem gives a characterization of open books supporting tight contact structures.

\begin{theorem}\label{thm:nrv-ot}
\cite[Theorem 1.1]{hkm1}
If $\phi$ is not right-veering then $(S,\phi)$ supports an overtwisted contact structure.
\end{theorem}

\begin{remark}
Conversely, Honda, Kazez and Mati\'c also prove that given an overtwisted contact structure $\xi$ there exists a non-right-veering open book $(S,\phi)$ supporting $\xi$ in \cite[p.444]{hkm1} where Eliashberg's classification of overtwisted contact structures \cite{el1} plays an important role. 
By Eliashberg's classification, an overtwisted contact structure admits an open book which is a negative stabilization of some open book.  Clearly such an open book has a non-right-veering monodromy. Therefore, a contact structure $\xi$ is tight if and only if every open book supporting $\xi$ is right-veering.
\end{remark}

\begin{proof}
If $\phi$ is not right-veering, then there exists a properly embedded oriented arc $\alpha \subset S$ such that $\phi(\alpha) > \alpha$. 
By \cite[Lemma 5.2]{hkm1}, there exists a sequence of properly embedded oriented arcs $\alpha_0, \ldots, \alpha_k$ such that 
\begin{itemize}
\item[(i)]
$\alpha_0,\dots,\alpha_k$ have the same initial point, $n \in \partial S$, 
\item[(ii)]
$\phi(\alpha)=\alpha_0 > \cdots > \alpha_k=\alpha$, 
\item[(iii)] 
consecutive $\alpha_i$ and $\alpha_{i+1}$ have disjoint interiors and distinct terminal points $p_i, p_{i+1}$.
\end{itemize}
Since $\phi(\alpha_k)=\alpha_0$ and $\phi= id$ near the bindings we have $p_0=p_k$. 
We may assume that: 
\begin{itemize}
\item[(iv)]
the terminal points $p_0, p_1, \ldots, p_{k-1} \in \partial S$ are mutually distinct.
\end{itemize}

Let $\beta_i$ (resp. $\check\beta_i$) be a sub-arc of $\alpha_i$ whose endpoints are $p_i$ and a point very close to $p_i$ (resp. $n$).
See Figure~\ref{arc-beta}.
\begin{figure}[htbp]
\begin{center}
\SetLabels
(0*.5) $n$\\
(1*.5) $p_i$\\
(.9*.1) $\beta_i$\\
(.6*.75) $\alpha_i$\\
(.5*.0) $\check{\beta_i}$\\
(-.05*0) $\partial S$\\
(1.05*.9) $\partial S$\\
\endSetLabels
\strut\AffixLabels{\includegraphics*[height=2cm]{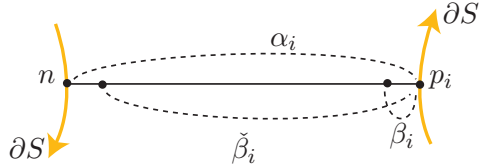}}
\caption{Arcs $\alpha_i, \beta_i,$ and $\check\beta_i$ oriented from $p_i$ to $n$.}\label{arc-beta}
\end{center}
\end{figure}
We orient $\alpha_i$ {\em against} the parametrization, i.e., the positive direction of $\alpha_i$ is from $p_i$ to $n$. 
The orientation of $\alpha_i$ induces those of $\beta_i$ and $\check\beta_i$.
We define sets of oriented arcs for $i=0, \ldots,k$: 
\[
A_i = \beta_0 \cup \cdots \cup \beta_{i-1} \cup \alpha_i \cup \beta_{i+1} \cup \cdots \cup \beta_k.
\]
Let $t_i = \frac{i}{k} \in [0,1]$. 
In the following, we construct an oriented surface $D_i$ properly embedded in the product region $S \times [t_i, t_{i+1}]$, where $i=0, \ldots, k-1$, such that $D_i \cap S_{t_i} =-A_i$ and $D_i \cap S_{t_{i+1}} =A_{i+1}$.

The surface $D_i$ consists of $k$ connected components; one non-product region and $k-1$ product regions defined by $\beta_j \times [t_i, t_{i+1}]$ where $j \neq i, i+1$ and $0\leq j\leq k$.

In the open book foliation of $D_i$ the point $n$ becomes a negative elliptic point, the points $p_0, \ldots, p_{k-1}$ become positive elliptic points, and the arc $\beta_j \times \{t\}$ becomes an a-arc in the page $S_t$. 
(See Prop. 2.2 of \cite{ik1-1} for the definition of a-arcs.) 
The non-product component of $D_i$ is defined by the movie presentation as sketched in Figure~\ref{hyp-point}. 
It is a saddle shape surface with a positive hyperbolic point $h_i$. 
\begin{figure}[htbp]
\begin{center}
\SetLabels
(.05*.95) (on page $S_{t_i})$\\
(-.01*.73) $n$\\
(.23*.9) $p_i$\\
(.28*.65) $p_{i+1}$\\
(.15*.8) $\alpha_i$\\
(.21*.7) $\beta_{i+1}$\\
(.36*.73) $n$\\
(.6*.9) $p_i$\\
(.65*.65) $p_{i+1}$\\
(.43*.79) $h_i$\\
(.8*.95) (on page $S_{t_{i+1}}$)\\
(.73*.73) $n$\\
(.98*.9) $p_i$\\
(1.02*.65) $p_{i+1}$\\
(.9*.7) $\alpha_{i+1}$\\
(.92*.82) $\beta_i$\\
(.2*.12) $n$\\
(.44*.31) $p_i$\\
(.49*.06) $p_{i+1}$\\
(.36*.2) $\alpha_i$\\
(.35*.0) $\check{\beta_i}$\\
(.27*.01) $+$\\
(.57*.12) $n$\\
(.81*.31) $p_i$\\
(.86*.06) $p_{i+1}$\\
(.72*.28) $\check{\beta_i}$\\
(.7*.05) $\alpha_{i+1}$\\
\endSetLabels
\strut\AffixLabels{\includegraphics*[height=60mm]{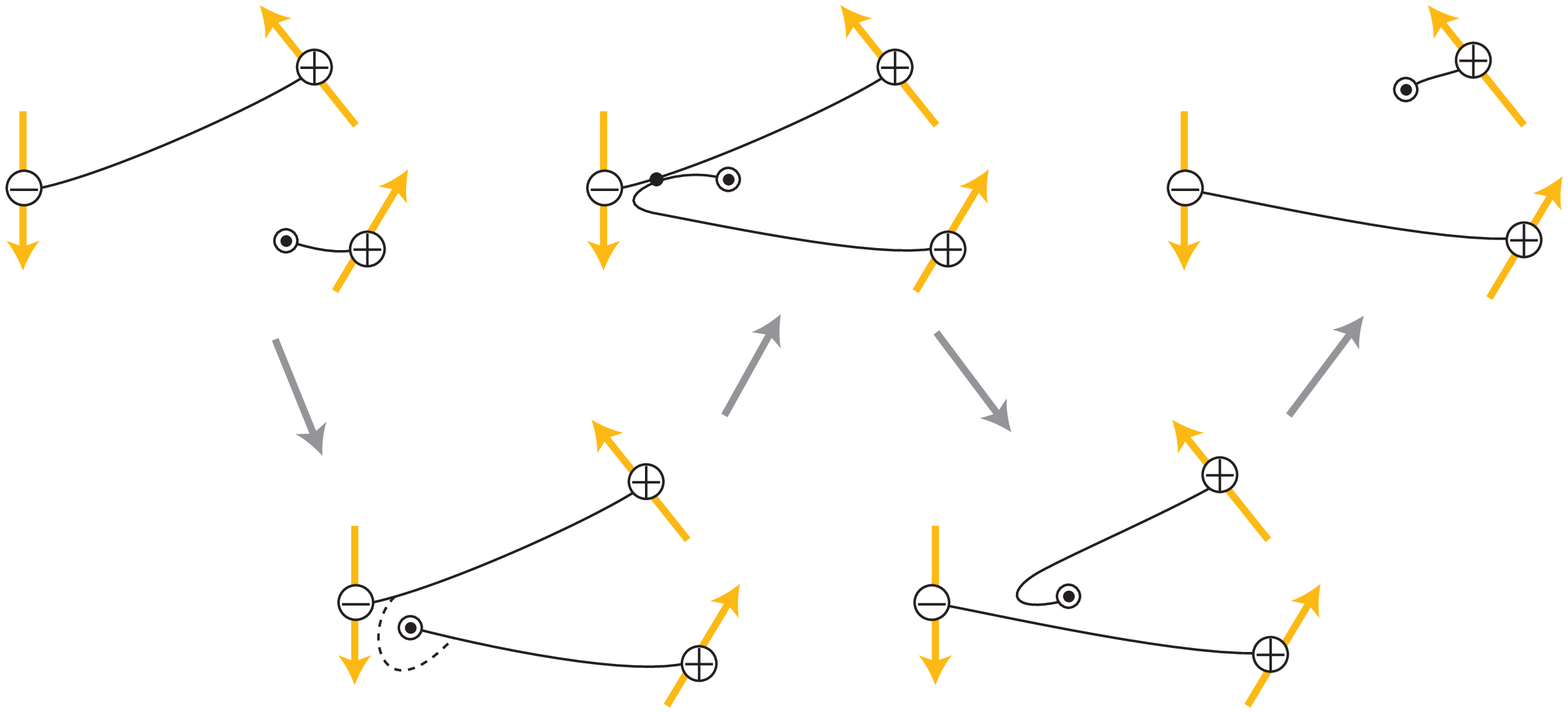}}
\caption{The non-product region of surface $D_i$ (movie presentation).}\label{hyp-point}
\end{center}
\end{figure}

Now we glue $D_i$ and $D_{i+1}$ along $A_{i+1} \subset S_{t_{i+1}}$ ($i=0, \ldots, k-2$) and obtain a surface $D_0 \cup \cdots \cup D_{k-1} \subset S \times [0,1]$ whose oriented boundary is $(-A_0) \cup A_k$.
Since arcs $\beta_1, \ldots, \beta_k$ are very close to $\partial S$ and $\phi=id$ near $\partial S$, we have $A_0 = \phi(A_k)$. 
So in the manifold $M_{(S, \phi)}$ we can identify $A_0$ and $A_k$ and obtain a surface which we denote by $D$.

The topological type of $D$ is the disc and its open book foliation $\F(D)$ is depicted in Figure~\ref{ot-disc}. Clearly our $D$ is a transverse overtwisted disc.
\begin{figure}[htbp]
 \begin{center}
 \SetLabels
(0.45*0.55) \large $n$\\
(0.8*0.25) $p_{0}=p_{k}$\\
(0.5*0.14) $h_{0}$\\
(0.38*0.11) $p_{1}$\\
(0.25*0.32) $h_{1}$\\
(0.2*0.65) $p_{2}$\\
(0.87*0.7) $p_{k-1}$\\
(0.8*0.43) $h_{k-1}$\\
\endSetLabels
\strut\AffixLabels{\includegraphics*[height=5cm]{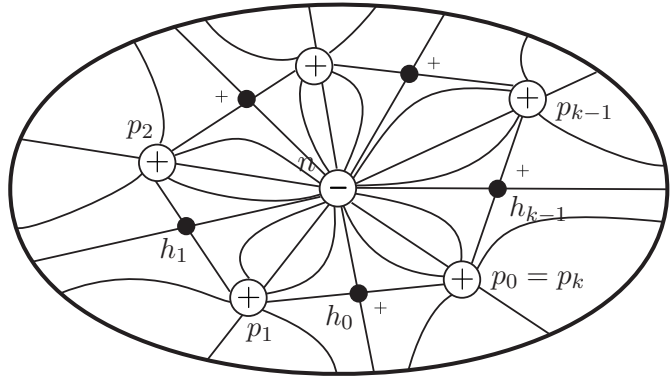}}
 \caption{The transverse overtwisted disc $D$.}\label{ot-disc}
  \end{center}
\end{figure}
\end{proof}

\begin{remark}
Honda, Kazez and Mati\'c's proof and our proof are based on the same combinatorial lemma \cite[Lemma 5.2]{hkm1} in order to use the assumption of right-veeringness. 
Our proof is more elementary and different from the original one that is written in the language of convex surface theory:

In \cite{hkm1}, they prove the {\em existence} of a bypass, half of an overtwisted disc, by applying \cite[Lemma 5.2]{hkm1}  and the right-to-life principal \cite[Lemma 2.9]{ho} \cite[Proposition 2.2]{hkm0} which involves the Legendrian realization principle and Eliashberg's classification of tight contact structures on  the 3-ball. 

On the other hand, we use \cite[Lemma 5.2]{hkm1} to explicitly {\em construct} a chain of positive elliptic and hyperbolic points surrounding the center negative elliptic point of an overtwisted disc. 
Hence our proof concretely visualizes an overtwisted disc. 
\end{remark}

\section{Overtwisted discs in right-veering open books}

{ \quad } \\ 
The converse of Theorem \ref{thm:nrv-ot} does not hold in general. 
Honda, Kazez and Mati\'c \cite{hkm1} show that {\em every} contact structure is supported by a right-veering open book:

Their argument is the following: 
Given a contact structure $(M, \xi)$ choose a compatible  open book $(S, \phi)$. 
For a boundary component $C$ of $S$ take two boundary-parallel arcs $a_{C}$ and $b_{C}$ such that the geometric intersection number $i(a_C, b_C)=2$. 
Apply positive stabilizations to $(S, \phi)$ along $a_{C}$ and $b_{C}$ for all the boundary components $C$ on which $\phi$ is non-right-veering. 
The new open book $(\widehat{S},\widehat{\phi})$ is now right-veering, see \cite[Prop.6.1]{hkm1}, and supports the same contact structure $\xi$.

Suppose that we start from a non-right-veering open book $(S,\phi)$, hence $\xi$ is overtwisted. 
In the following we concretely describe how is an overtwisted disc embedded with respect to the stabilized right-veering open book $(\widehat{S},\widehat{\phi})$. 

\begin{example}[Overtwisted discs in Honda-Kazez-Mati\'c's stabilizations]

Let $(S,\phi)$ be a not-right-veering open book. 
By the proof of Theorem \ref{thm:nrv-ot} one can construct an overtwisted disc, $D$, in $M_{(S,\phi)}$.
The open book foliation of $D$ has a unique negative elliptic point, say $n$, that lies on the binding component $C \subset \partial S$,  
and $k$ positive elliptic points $p_{1},\ldots,p_{k}$ and $k$ positive hyperbolic points $h_{1},\ldots,h_{k}$.
Let $S_{t_{i}}$ be the singular fiber that contains $h_{i}$.  We may assume 
\[ 0<t_1 <t_2 <\cdots < t_k < \frac{1}{2} <1.\]
For $t \in [0,1)$ let $b_{t} \in S_{t}\cap D$ be the b-arc (cf. Prop.2.2 of \cite{ik1-1}) that ends at the point $n$.

Now we apply Honda-Kazez-Mati\'c's stabilization to get a right-veering open book $(\widehat{S},\widehat{\phi})$. The monodromy $\widehat{\phi}$ satisfies $\widehat{\phi}= 
T_{\beta}\circ T_{\alpha} \circ \phi$, 
where $\alpha$ and $\beta$ are  core circles of the annuli plumbed along $a_C$ and $b_{C}$ and $T_\alpha, T_\beta$ are positive Dehn twists along $\alpha, \beta$.

We will construct an overtwisted disc $\widehat D$ by giving a movie presentation relative to the open book $(\widehat{S},\widehat{\phi})$. 
For sake of simplicity we assume that: 
\begin{itemize}
\item 
$\phi$ is non-right-veering {\em only} along $C$. 
\item
$p_k \in \partial S \setminus C$. 
\end{itemize}
In the general case a construction of $\widehat D$ is similar but more complicated. It is obtained as an application of arguments in \cite{ik3}.

Choose stabilization arcs $a_{C}$ and $b_{C}$ such that $i(\alpha, b_{1/2}) =1$ and $i(\beta, b_{1/2}) =0$ as shown in Figure \ref{fig:hkmexam1} (a).  
Such arcs can always be found by the assumptions  above. 
\begin{figure}[htbp]
\begin{center}
 \SetLabels
(0.05*0.95) (a)\\
(0.05*0.63) (b)\\
(0.55*0.63) (c)\\
(0.05*0.28) (d)\\
(0.11*0.81) $n$\\
(.6*.86) new b-arc\\
(.15*.86) $\widehat b_t$\\
(.15*.25) $\widehat b_1$\\
(0.18*0.94) $\beta$ \\
(0.25*0.92) $\alpha$ \\
(0.49*0.82) $p'$ \\
(0.38*0.86) $n'$ \\
(0.32*0.56) $h_{-}$ \\
(0.83*0.44) $h_{+}$ \\
(0.48*0.71) $S_{t\in [0, 1/2]}$\\
(0.48*0.02) $S_{1}$\\
(.25*.7) $C$\\
\endSetLabels
\strut\AffixLabels{\includegraphics*[width=100mm]{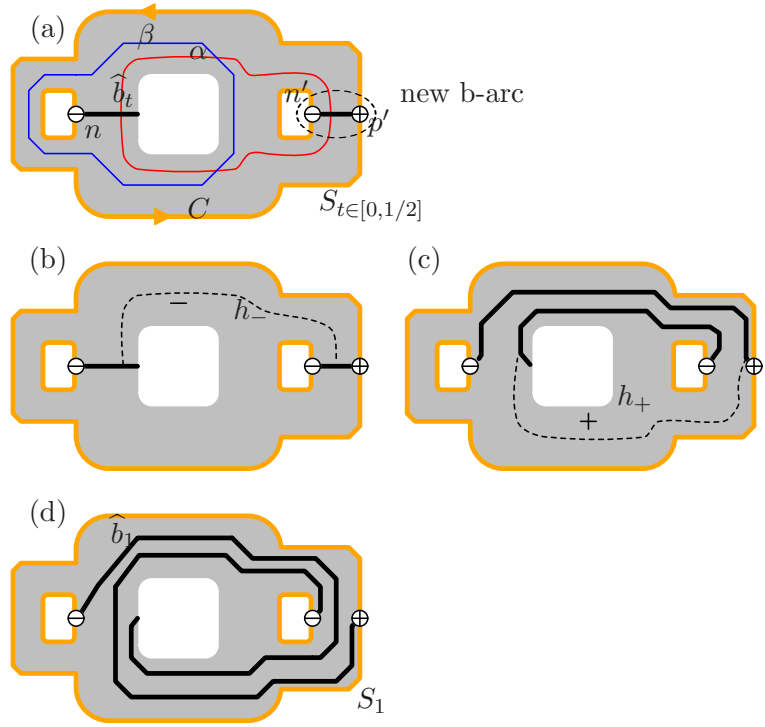}}
\caption{A movie presentation of the overtwisted disc $\widehat D$. The shaded region is the union of a collar neighborhood of $C$ and the two plumbed annuli. All the a-arcs are omitted for simplicity. } 
\label{fig:hkmexam1}
\end{center}
\end{figure}

To the region $\{ \widehat S_t \ | \ t \in [0,\frac{1}{2}]\}$ add a continuous family of b-arcs that are co-cores of the annuli plumbed along $a_C$, see Figure \ref{fig:hkmexam1} (a). We denote the positive and negative elliptic points which are the endpoints of the newly added $b$-arcs by $p'$ and $n'$. Except this family of b-arcs, the movie presentation of $\widehat D$ in the interval $[0,\frac{1}{2}]$ is the same as that of $D$. 
Let $\widehat b_t := b_t$.

In the interval $[\frac{1}{2},1)$, $\widehat{D}$ is described as in the passage of (b)$\to$(c)$\to$(d) of Figure \ref{fig:hkmexam1}. 
We form one negative and one positive hyperbolic points $h_{-}$ and $h_{+}$ as in (b) and (c), respectively. The describing arcs of $h_{-}$ and $h_{+}$ are parallel to $\alpha$.
We have $\widehat b_1 = T_\alpha^{-1} (b_1)$. 
Note that 
$\widehat \phi(\widehat b_1) 
= (T_{\beta} \circ T_{\alpha} \circ \phi) (\widehat b_1) 
= T_{\beta} \circ T_{\alpha} (T_\alpha^{-1} b_0)
=T_{\beta}  (b_0) =  b_0 
= \widehat b_0$, 
moreover $\widehat \phi(\widehat D \cap S_1)=\widehat D \cap S_0$ 
so this movie presentation indeed defines an overtwisted disc in $M_{(\widehat{S},\widehat{\phi})}$.

The open book foliations of the overtwisted discs $D$ and $\widehat{D}$ are depicted in Figure \ref{fig:hkmotdisc1}, where $\F(\widehat{D})$ is obtained by inserting two bb-tiles of opposite signs into the shaded region of $\F(D)$ that is bounded by $b_{\frac{1}{2}}$ and $b_{1}$. 
\begin{figure}[htbp]
\begin{center}
\SetLabels
(.27*.57) $b_{\frac{1}{2}}$\\
(.27*.42) $b_1$\\
(.84*.7) $b_{\frac{1}{2}}$\\
(.83*.3) $b_1$\\
(0.16*0.56) $n$\\
(0.57*0.56) $n$\\
(0.68*0.54) $p'$\\
(0.86*0.54) $n'$\\
(0.76*0.33) $h_{+}$\\
(0.76*0.63) $h_{-}$\\
(0.33*0.38) \small $p_{k}$\\
(0.94*0.38) \small $p_{k}$\\
(0.18*0.21) \small $h_{2}$\\
(0.25*0.13) \small $p_{1}$\\
(0.25*0.86) \small $p_{k-1}$\\
(0.16*0.76) \small $h_{k-1}$\\
(0.31*0.72) \small $h_{k}$\\
\endSetLabels
\strut\AffixLabels{\includegraphics*[width=125mm]{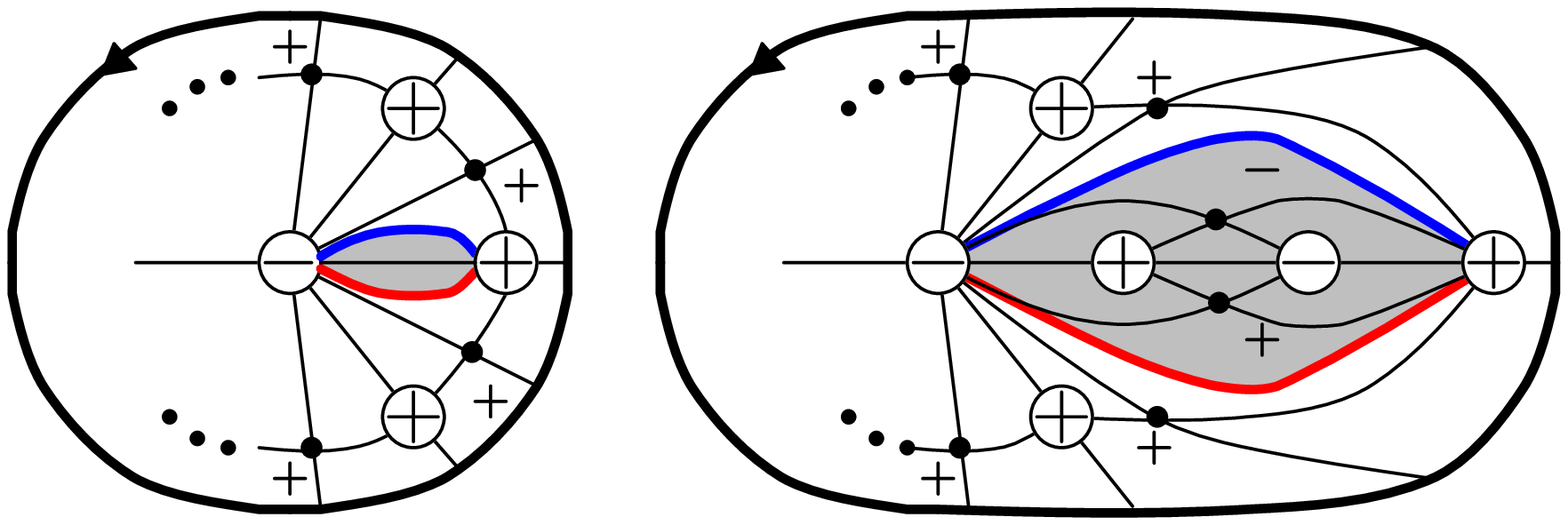}}
\caption{$\F(D)$ and $\F(\widehat{D})$}
\label{fig:hkmotdisc1}
\end{center}
\end{figure}

Strictly speaking, the disc $\widehat{D}$ is not a transverse overtwisted disc since the condition (3) of Definition \ref{def:trans-ot-disc}  is not satisfied. 
However, applying the same technique we use in the alternative proof in \cite{ik1-1} of the Bennequin-Eliashberg inequality \cite{Ben, el2}, the condition (3) will be satisfied. Thus we can regard $\widehat{D}$ as an overtwisted disc. 
\end{example}

\section{Generalization of Lekili and Lisca's examples}

In \cite[Question 6.2]{hkm1} Honda, Kazez and Mati\'c ask whether a right-veering and non-destabilizable open book always supports a tight contact structure.
Lekili \cite{le} and Lisca \cite{lis} negatively answer the question by constructing examples. They study open book decompositions of $3$-manifolds whose tight contact structures are well-studied and classified 
(in \cite{le} Poincar\'e homology $3$-spheres, and in \cite{lis} lens spaces). 
In both constructions the most technical points are showing that their open books indeed support overtwisted contact structures. 
Advanced tools such as Ozsv\'ath-Sz\'abo's Heegaard Floer invariants and properties of planar open books enable them to overcome the difficulty. 

We generalize Lekili and Lisca's examples in Theorem~\ref{theorem:rvot} below. 
Our proof of overtwistedness is direct and does not require any knowledge of classification of tight contact structures of ambient manifolds or Ozsv\'ath-Sz\'abo's invariants.

\begin{theorem}
\label{theorem:rvot}
Let $S$ be a $2$-sphere with four holes. Let $a,b,c,d,e$ be simple closed curves on $S$ as shown in Figure~\ref{4psphere}.
Let 
$\Phi_{h,i,k} = T_{a}^{h}T_{b}^{i}T_{c}T_{d}T_{e}^{-k-1}$
where $T_{x}$ $(x=a,b,c,d,e)$ denotes the right-handed Dehn twist along $x$.
Then for all $h,i,k \geq 1$, $\Phi_{h,i,k}$ is right-veering and the open book $(S,\Phi_{h,i,k})$ is non-destabilizable and supports an overtwisted contact structure.
\begin{figure}[htbp]
\begin{center}
\strut\AffixLabels{\includegraphics*[height=35mm]{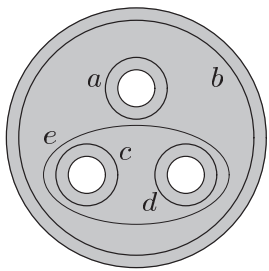}}
\caption{The surface $S$.}\label{4psphere}
\end{center}
\end{figure}
\end{theorem}

\begin{remark}
Lekili's examples \cite[Theorem 1.2 and Remark 4.1]{le}  are $\Phi_{2,i,1}$ ($i \leq 5$), and Lisca's examples \cite[Theorem 1.1]{lis} are $\Phi_{h,1,l}$ $(h,l >0)$. 
\end{remark}

\begin{figure}[htbp]
\begin{minipage}{0.5\textwidth}
\SetLabels
(-0.03*.95)   $(7)$\\
(.40*.50) {\footnotesize $2k+1$}\\
(.51*.55)   \tiny $1$\\
(.51*.50)   \tiny $1$\\
(.56*.44) \tiny $2k-1$\\
(.25*.45) \tiny $1$\\
(.25*.5)  \tiny $1$\\
(.22*.56)  \tiny $2k-1$\\
(.32*.72)  $h$\\
(.21*.8)   \tiny $h-1$\\
(.33*.8)    \tiny $1$\\
(.26*.65)   \tiny $1$\\
(.39*.65)    \tiny $h-1$\\
(.4*.06)  $i$\\
(.52*.04)  \tiny $1$\\
(.55*.09)  \tiny $i-1$\\
(.3*.09)    \tiny $1$\\
(.26*.04)  \tiny $i-1$\\
(1.12*.05)   page $S_{1}$\\
\endSetLabels
\strut\AffixLabels{\includegraphics*[scale=0.5, width=54mm]{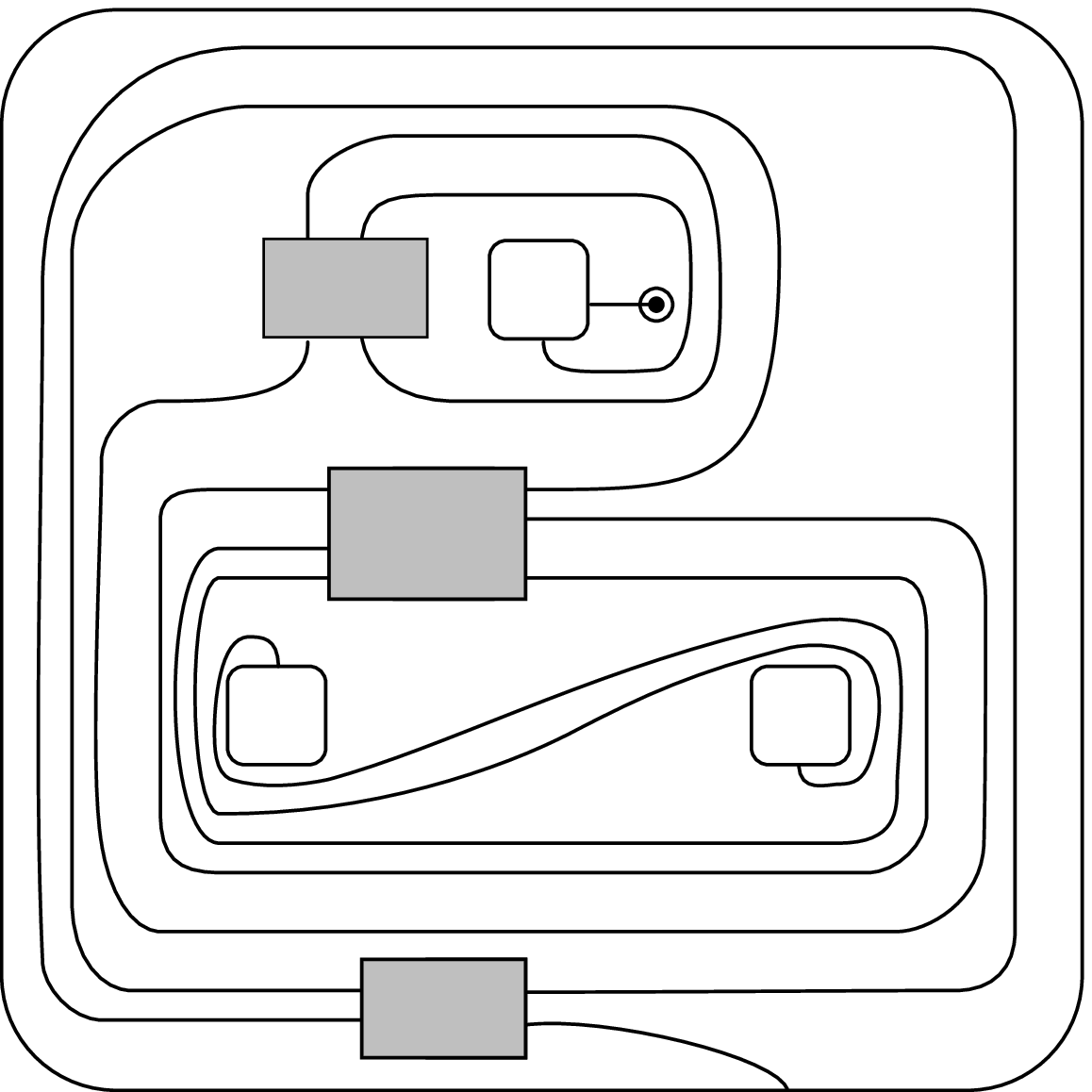}}
 \end{minipage}
  \begin{minipage}{0.48\textwidth}
\begin{itemize}
\item The leaves in the page $S_{1}$. It satisfies $\Phi_{h,i,k}(D \cap S_1) = D \cap S_0$.
\end{itemize}
 \end{minipage}
 \caption{Movie presentation of a transverse overtwisted disc in $(S,\Phi_{h,i,k})$.}\label{fig:Liscaexample}
\end{figure}

\begin{figure}[htbp]
\begin{center}
\SetLabels 
\endSetLabels
\strut\AffixLabels{\includegraphics*[height=45mm]{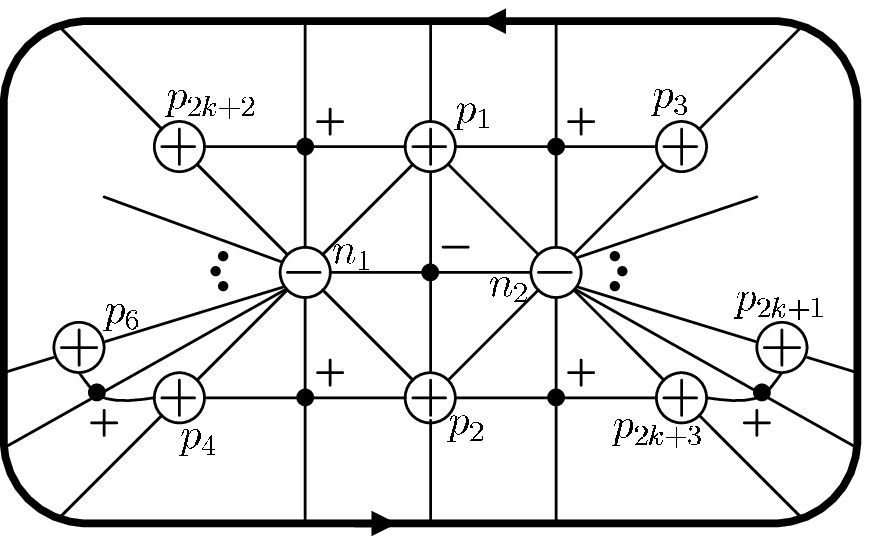}}
 \end{center}
 \caption{The open book foliation $\F(D)$.}
 \label{fig:Liscafoliation}
 \end{figure}

\begin{proof}
We owe \cite{le} \cite{lis} the proof that $\Phi_{h,i,k}$ is right-veering and non-destabilizable. 
Hence we only show that $(S,\Phi_{h,i,k})$ supports an overtwisted contact structure.
We define a transverse  overtwisted disc $D$ in the open book $(S,\Phi_{h,i,k})$ by the movie presentation as in Sketches (1)--(7) of Figure \ref{fig:Liscaexample}. 
For example, Sketch (1) depicts the page $S_0$ with the set of arcs $D\cap S_0$. 
On $\partial S$ there are two negative elliptic points $n_{1}$, $n_{2}$,  and $(2k+3)$ positive elliptic points $p_{1},\ldots,p_{2k+3}$. 
The movie presentation shows that $\F(D)$ contains two negative hyperbolic points and $2k+3$ positive hyperbolic points. 
Note that $\Phi_{h,i,k}(D \cap S_1)=D \cap S_0$. 

The corresponding open book foliation $\F(D)$ is depicted in Figure~\ref{fig:Liscafoliation}. 
We can verify that $D$ meets the conditions in Definition~\ref{def:trans-ot-disc} of transverse overtwisted discs. 
\end{proof}

\begin{figure}[htbp]
\begin{minipage}{0.5\textwidth}
\SetLabels
(-0.03*.95)    $(1)$\\
(.38*.7)    $p_{1}$\\
(.67*-.03)  $p_{2}$\\
(.55*.63)   $p_{2k+3}$\\
(.95*.52)  $p_{3}$\\
(.99*.7) $p_{2k+2}$\\
(.25*.36)     $n_{1}$\\
(.67*.33)   $n_{2}$\\
(1.00*.05) page $S_{0}$\\
(.45*.27) ${\color{red} -}$\\
\endSetLabels
\strut\AffixLabels{\includegraphics*[scale=0.5, width=63mm]{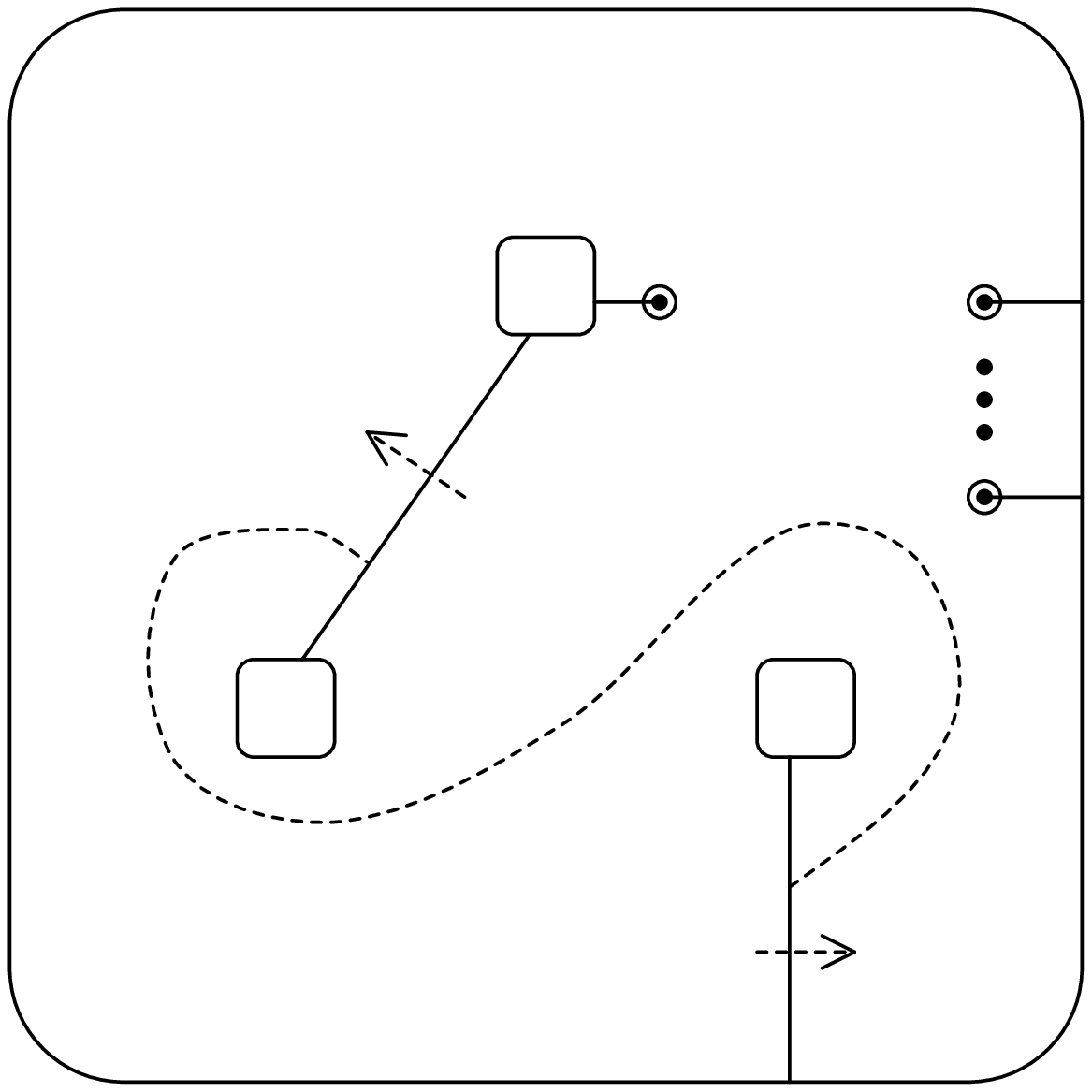}}
 \end{minipage}
 \begin{minipage}{0.48\textwidth}
\begin{itemize}
\item Leaves in the fiber $S_{0}$ consist of two $b$-arcs and $(2k+1)$ $a$-arcs.
\item The dashed describing arc corresponds to a negative hyperbolic point.
\end{itemize}
\end{minipage}
\end{figure}

\begin{figure}[htbp]
\begin{minipage}{0.5\textwidth}
\noindent
\SetLabels
(-0.03*.95)   $(2)$\\
(.38*.7)    $p_{1}$\\
(1.07*.54)  $p_{3}$ (or $p_{2i+1})$\\
(.67*-.03)  $p_{2}$\\
(.25*.38)     $n_{1}$\\
(.67*.34)   $n_{2}$\\
(.45*.55) ${\color{red} +}$\\
\endSetLabels
\strut\AffixLabels{\includegraphics*[scale=0.5, width=63mm]{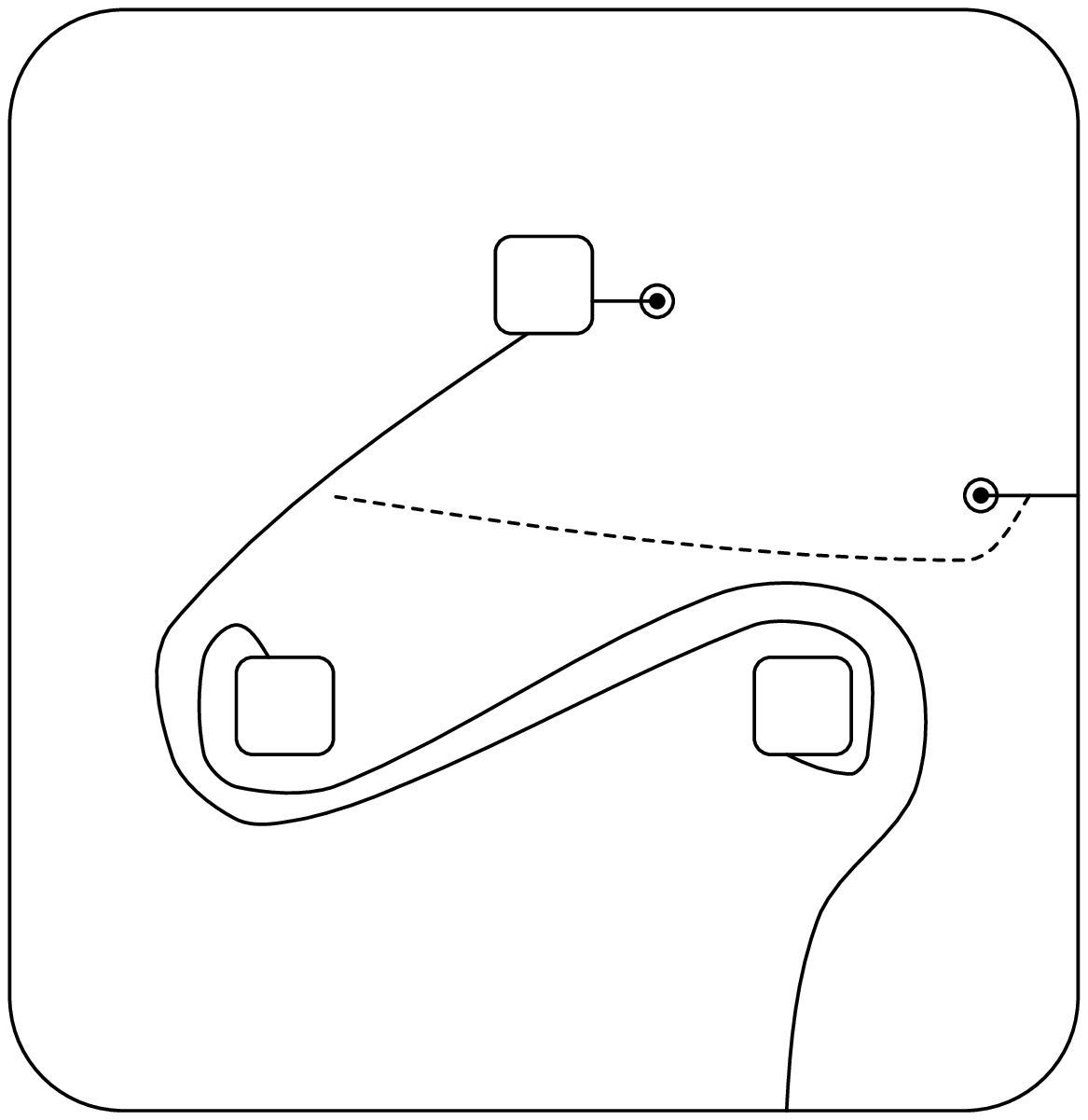}}
 \end{minipage}
  \begin{minipage}{0.48\textwidth}
  \begin{itemize}
\item In Sketches $(2),\dots,$(7), we omit some $a$-arcs emanating from $p_{3},\ldots,p_{2k+2}$ if they are not involved in producing hyperbolic singularities. 
\item The describing arc joining the $b$-arc emanating from $n_2$ and the $a$-arc from $p_3$ (or $p_{2i+1}$ in the $i$-th iteration) represents  a positive hyperbolic point. 
\end{itemize}
 \end{minipage}
\end{figure}

\begin{figure}[htbp]
\begin{minipage}{0.5\textwidth}
\noindent
\SetLabels
(-0.03*.95)   $(3)$\\
(1.1*.54)  $p_3$ (or $p_{2i+1})$\\
(1.1*.62)   $p_{4}$ (or $p_{2i+2})$\\
(.25*.38)     $n_{1}$\\
(.67*.34)   $n_{2}$\\
(.7*-.03)  $p_{2}$\\
(.35*.6) ${\color{red} +}$\\
\endSetLabels
\strut\AffixLabels{\includegraphics*[scale=0.5, width=63mm]{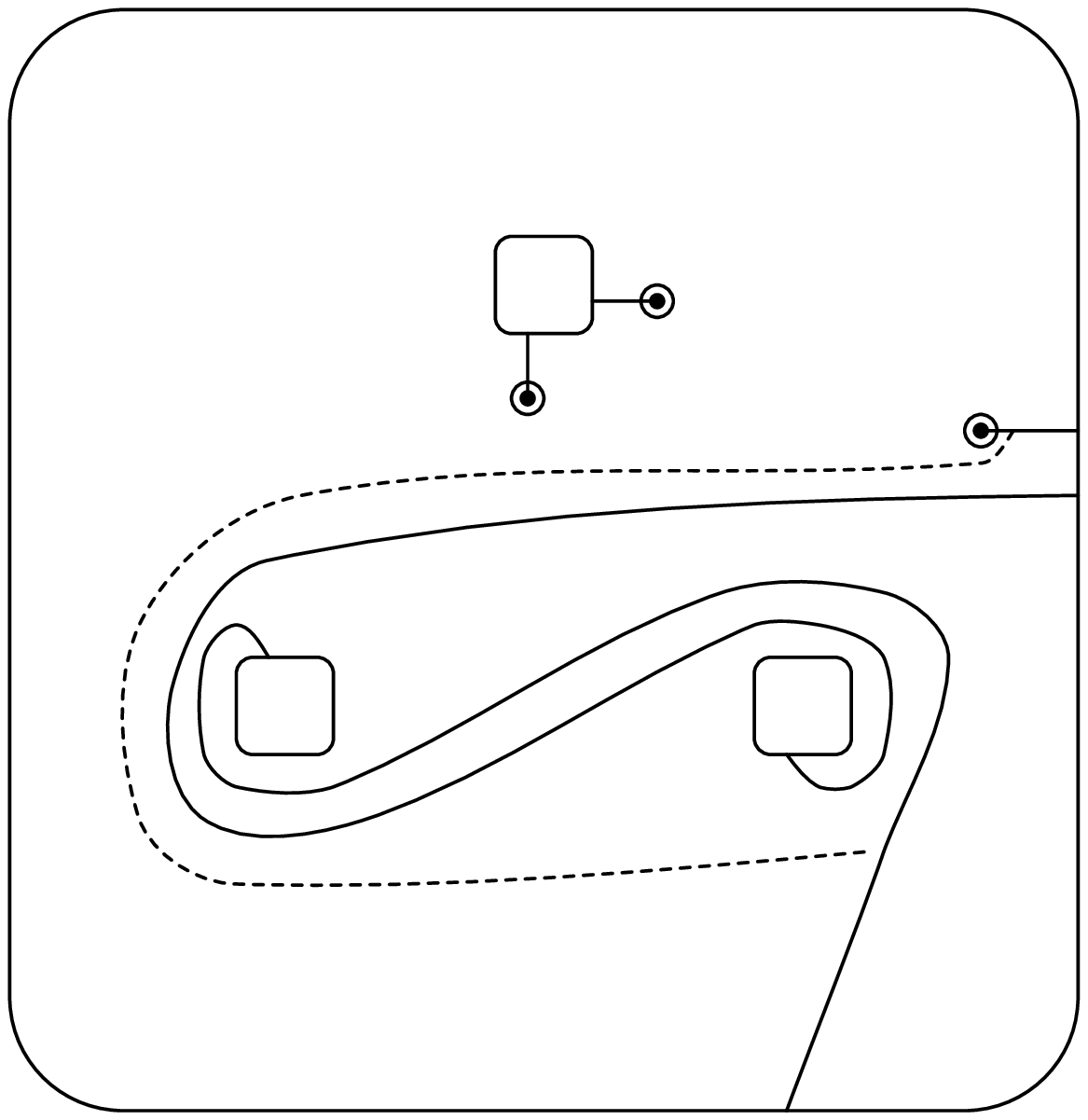}}
 \end{minipage}
  \begin{minipage}{0.48\textwidth}
\begin{itemize}
\item The describing arc joining the $b$-arc emanating from $n_{1}$ and the $a$-arc from $p_{4}$ (or, $p_{2i+2}$ in the $i$-th iteration) represents a positive hyperbolic point.
\item Iterate the Steps (2) and (3) for $k$ times.
\end{itemize}
 \end{minipage}
\end{figure}

\begin{figure}[htbp]
\begin{minipage}{0.5\textwidth}
\SetLabels
(-0.03*.95)   $(4)$\\
(.4*.1) ${\color{red} +}$\\
(1.1*.55)  $p_{2k+1}$\\
(.5*.55) \tiny $1$\\
(.5*.5)   \tiny $1$\\
(.55*.45) \tiny $2k-2$\\
(.25*.45) \tiny $1$\\
(.25*.5)  \tiny $1$\\
(.2*.55)  \tiny $2k-2$\\
(.25*.78) \tiny $h-1$\\
(.35*.78) \tiny $1$\\
(.27*.66) \tiny $1$\\
(.37*.66) \tiny $h-1$\\
(1.1*.65)   $p_{2k+2}$\\
(.6*.78)     $p_{2k+3}$\\
(.4*.50) $2k$\\
(.32*.72)   $h$\\
(.75*.32)   $n_{2}$\\
\endSetLabels
\strut\AffixLabels{\includegraphics*[scale=0.5, width=55mm]{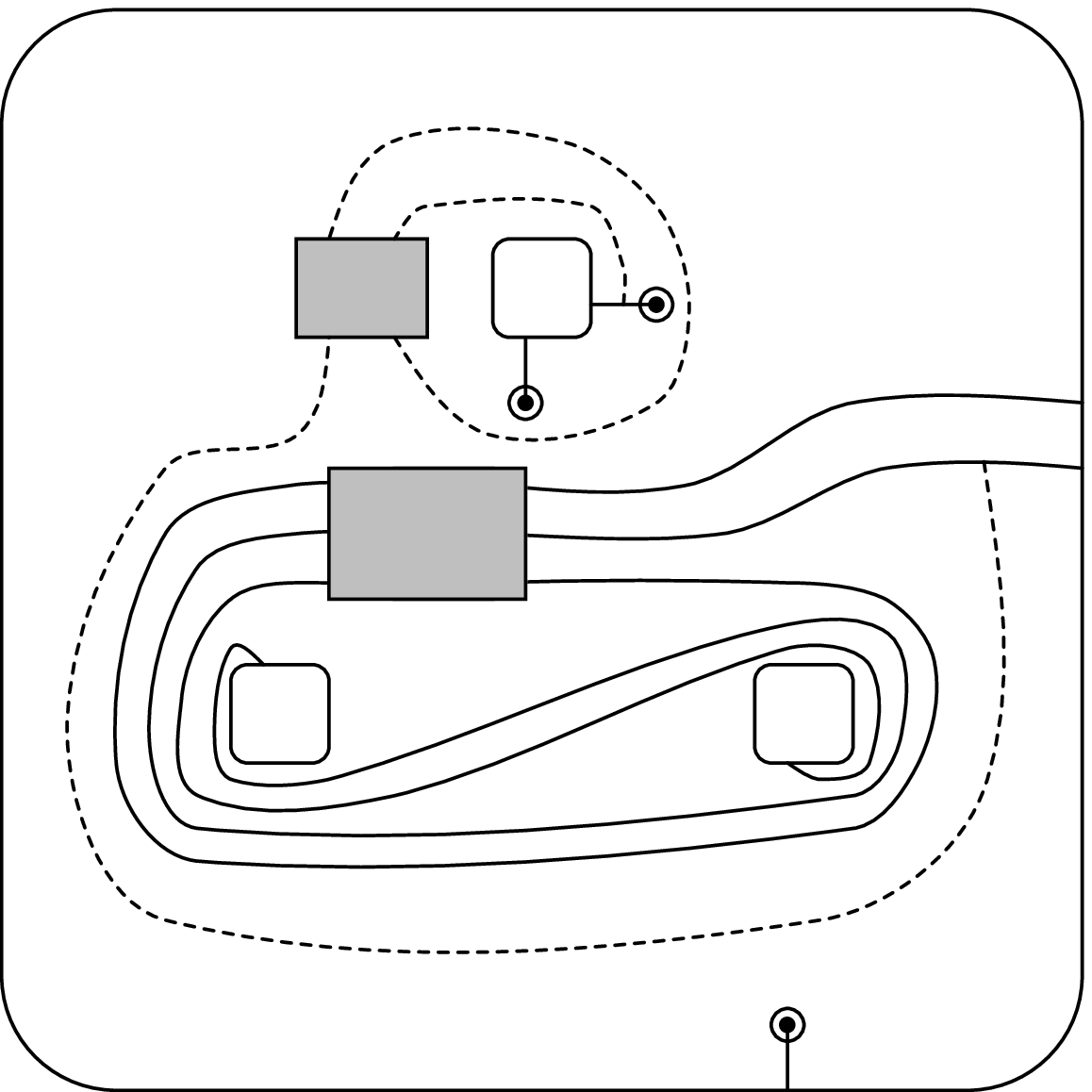}}
 \end{minipage}
  \begin{minipage}{0.48\textwidth}
\begin{itemize}
\item The shaded boxes labeled $h, 2k$ contain parallel $h, 2k$ arcs. 
\item The edges out of the shaded boxes are weighted as indicated. 
\item The describing arc joins the $b$-arc from $n_{2}$ to $p_{2k+1}$ and the $a$-arc emanating from $p_{2k+3}$ and it represents a positive hyperbolic point. 
\end{itemize}
 \end{minipage}
\end{figure}

\begin{figure}[htbp]
\begin{minipage}{0.5\textwidth}
\SetLabels
(-0.03*.95)  $(5)$\\
(.4*.09)       ${\color{red} +}$\\
(.5*.72)       $p_1$\\
(1.1*.65)     $p_{2k+2}$\\
(.39*.50)  $2k$\\
(.23*.55) \tiny $2k-2$\\
(.29*.5) \tiny $1$\\
(.29*.45) \tiny $1$\\
(.5*.55) \tiny $1$\\
(.5*.5) \tiny $1$\\
(.55*.45) \tiny $2k-1$\\
(.3*.73)     $2h$\\
(.26*.66) \tiny $1$\\
(.31*.66) \tiny $1$\\
(.2*.8) \tiny $h-1$\\
(.45*.88) \tiny $h-1$\\
(.35*.8) \tiny $1$\\
(.39*.8) \tiny $1$\\
(.26*.36)     $n_{1}$\\
\endSetLabels
\strut\AffixLabels{\includegraphics*[scale=0.5, width=55mm]{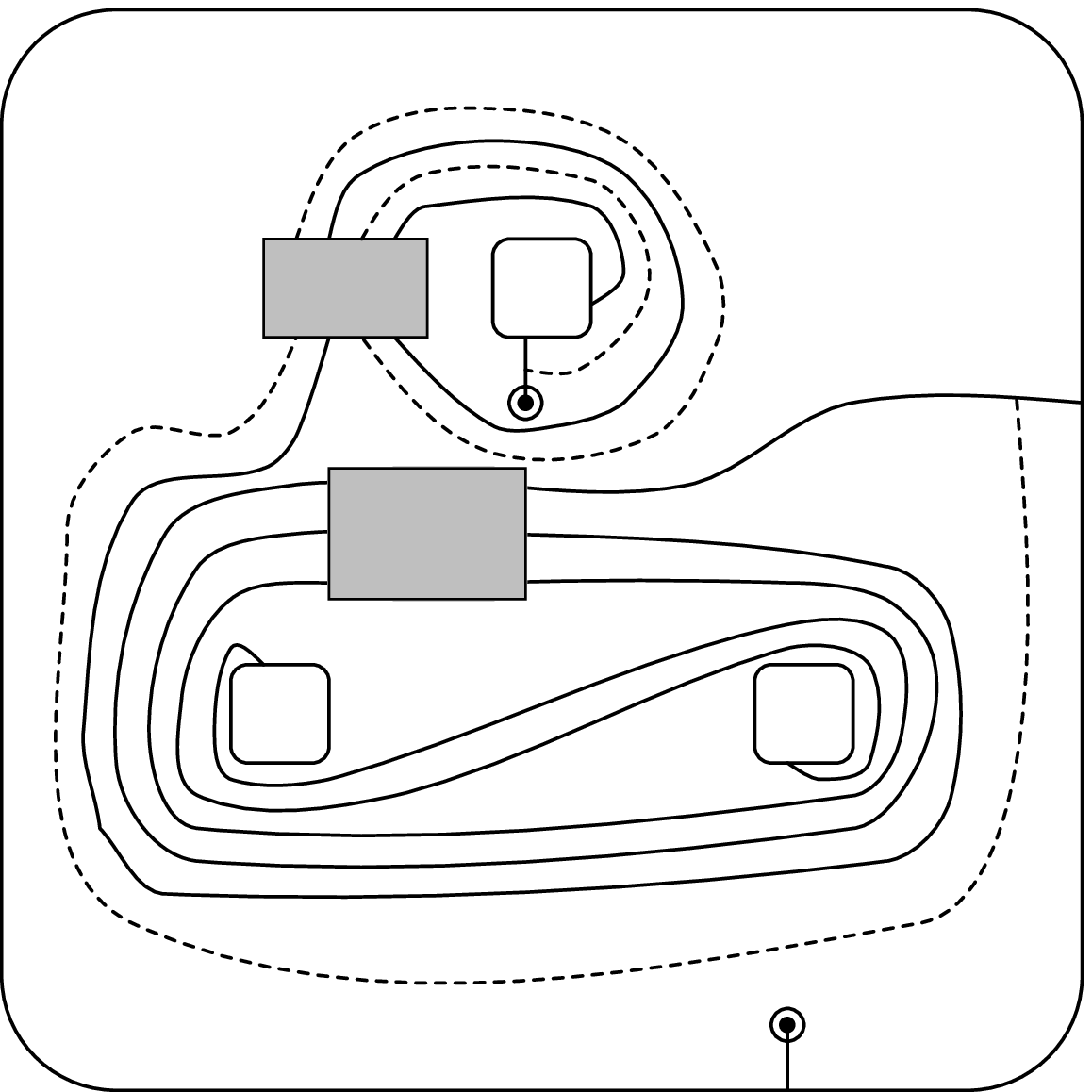}}
 \end{minipage}
  \begin{minipage}{0.48\textwidth}
\begin{itemize}
\item The describing arc joins the $b$-arc from $n_{1}$ to $p_{2k+2}$ and the $a$-arc emanating from $p_{1}$ and it represents a positive hyperbolic point. 
\end{itemize}
 \end{minipage}
\end{figure}

\begin{figure}[htbp]
\begin{minipage}{0.5\textwidth}
\SetLabels
(-0.03*.95)   $(6)$\\
(.36*.50)  $2k$\\
(.51*.51)   \tiny $2$\\
(.56*.44) \tiny $2k-2$\\
(.25*.45) \tiny $1$\\
(.25*.5)  \tiny $1$\\
(.2*.54)  \tiny $2k-2$\\
(.32*.71)      $2h$\\
(.21*.8)   \tiny $2h-2$\\
(.33*.8)    \tiny $1$\\
(.38*.8)    \tiny $1$\\
(.26*.65)   \tiny $2$\\
(.39*.65)    \tiny $2h-2$\\
(.39*.06)  $i$\\
(.52*.04)  \tiny $1$\\
(.55*.09)  \tiny $i-1$\\
(.3*.09)    \tiny $1$\\
(.26*.04)  \tiny $i-1$\\
(.65*.78)     $p_{2k+3}$\\
(.75*-.05)  $p_{2}$\\
(.75*.32)   $n_{2}$\\
(.9*.56)       ${\color{red} +}$\\
\endSetLabels
\strut\AffixLabels{\includegraphics*[scale=0.5, width=55mm]{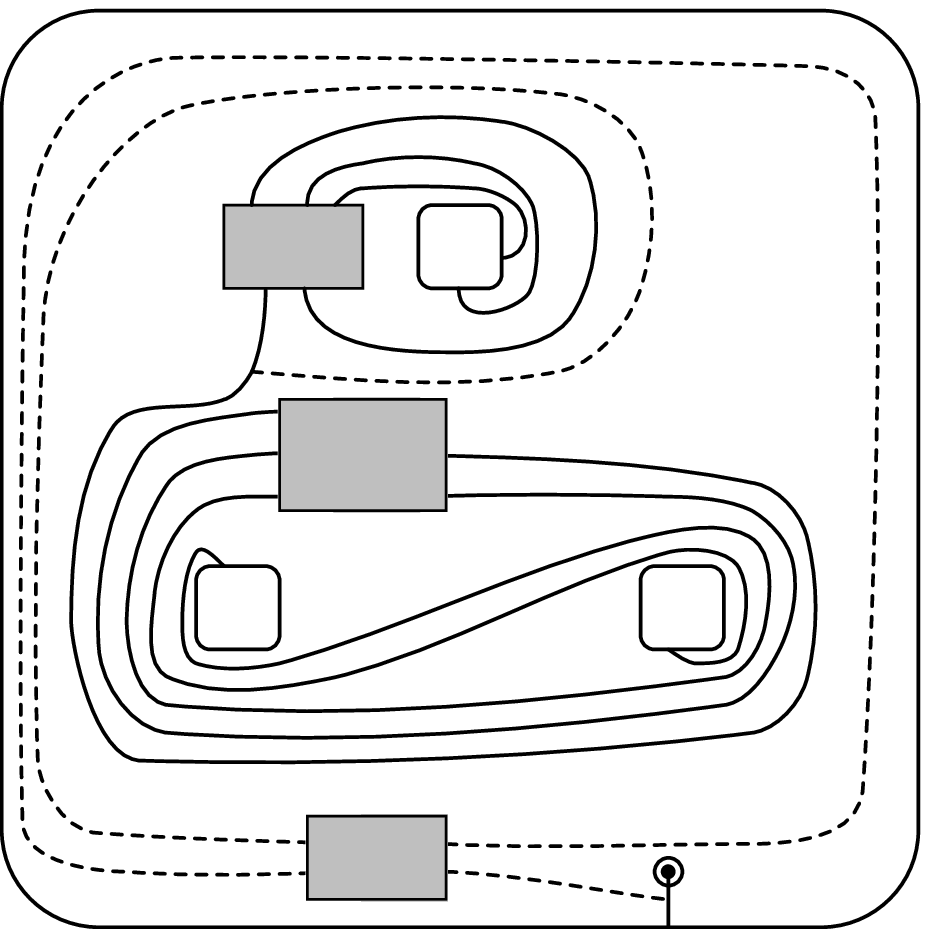}}
 \end{minipage}
  \begin{minipage}{0.48\textwidth}
\begin{itemize}
\item The describing arc joins the $b$-arc from $n_2$ to $p_{2k+3}$ and the $a$-arc emanating from $p_{2}$ and it represents a positive hyperbolic point. 
\end{itemize}
 \end{minipage}
\end{figure}

\begin{remark}
After the announcement of this example (on December 26, 2011), Kazez and Roberts \cite{KR} found more counterexamples to the conjecture of Honda, Kazez and Mati\'c. 
\end{remark}

\section*{Acknowledgement}
The first author was supported by JSPS Research Fellowships for Young Scientists.  
The second author was partially supported by NSF grants DMS-0806492 and DMS-1206770.

\end{document}